\documentclass[a4paper, 11pt,twoside]{article}
\usepackage{palatino}
\usepackage{textcomp}
\usepackage{amssymb}
\usepackage{graphicx}
\usepackage{epsfig}
\usepackage{epstopdf}
\usepackage{tikz}
\usepackage{mathrsfs}
\usepackage{picins}
\usepackage{mathdots}
\usepackage[english]{babel}

\newtheorem{theorem}{Theorem}
\newtheorem{proposition}{Proposition}
\newtheorem{lemma}{Lemma}
\newtheorem{corollary}{Corollary}
\newtheorem{definition}{Definition}
\newtheorem{example}{Example}

\newtheorem{property}{Property}

\newcommand{\subplus}{\begin{tikzpicture}[scale=0.08]\draw (0,0)--(2,0)--(2,2)--(0,2)--(0,0);\draw (1,0)--(1,2); \draw (0,1)--(2,1);\end{tikzpicture}}
\newcommand{\plus}{\begin{tikzpicture}[scale=0.11]\draw (0,0)--(2,0)--(2,2)--(0,2)--(0,0);\draw (1,0)--(1,2); \draw (0,1)--(2,1);\end{tikzpicture}}
\newcommand{\threetop}{\begin{tikzpicture}[scale = 0.11]\draw (0,0)--(2,0)--(2,2)--(0,2)--(0,0); \draw (0,1)--(2,1);\draw (1,0)--(1,1);\end{tikzpicture}}
\newcommand{\threeright}{\begin{tikzpicture}[scale = 0.11]\draw (0,0)--(2,0)--(2,2)--(0,2)--(0,0);\draw (1,0)--(1,2);\draw (0,1)--(1,1);\end{tikzpicture}}
\newcommand{\threeleft}{\begin{tikzpicture}[scale = 0.11]\draw (0,0)--(2,0)--(2,2)--(0,2)--(0,0);\draw (1,0)--(1,2);\draw (2,1)--(1,1);\end{tikzpicture}}
\newcommand{\threebottom}{\begin{tikzpicture}[scale = 0.11]\draw (0,0)--(2,0)--(2,2)--(0,2)--(0,0);\draw (0,1)--(2,1);\draw (1,2)--(1,1);\end{tikzpicture}}
\newcommand{\vertical}{\begin{tikzpicture}[scale = 0.11]\draw (0,0)--(2,0)--(2,2)--(0,2)--(0,0);\draw (1,0)--(1,2);\end{tikzpicture}}
\newcommand{\horizontal}{\begin{tikzpicture}[scale = 0.11]\draw (0,0)--(2,0)--(2,2)--(0,2)--(0,0);\draw (0,1)--(2,1);\end{tikzpicture}}
\newcommand{\tlandbr}{\begin{tikzpicture}[scale = 0.11]\draw (0,0)--(2,0)--(2,2)--(0,2)--(0,0);\draw (0,1) arc (270:360:1);\draw (2,1) arc (90:180:1); \end{tikzpicture}}
\newcommand{\trandbl}{\begin{tikzpicture}[scale = 0.11]\draw (0,0)--(2,0)--(2,2)--(0,2)--(0,0);\draw (1,2) arc (180:270:1);\draw (1,0) arc (0:90:1);\end{tikzpicture}}
\newcommand{\topleft}{\begin{tikzpicture}[scale = 0.11]\draw (0,0)--(2,0)--(2,2)--(0,2)--(0,0);\draw (0,1) arc (270:360:1); \end{tikzpicture}}
\newcommand{\topright}{\begin{tikzpicture}[scale = 0.11]\draw (0,0)--(2,0)--(2,2)--(0,2)--(0,0);\draw (1,2) arc (180:270:1); \end{tikzpicture}}
\newcommand{\bottomleft}{\begin{tikzpicture}[scale = 0.11]\draw (0,0)--(2,0)--(2,2)--(0,2)--(0,0);\draw (1,0) arc (0:90:1); \end{tikzpicture}}
\newcommand{\bottomright}{\begin{tikzpicture}[scale = 0.11]\draw (0,0)--(2,0)--(2,2)--(0,2)--(0,0);\draw (2,1) arc (90:180:1); \end{tikzpicture}}
\newcommand{\disconnected}{\begin{tikzpicture}[scale = 0.11]\draw (0,0)--(2,0)--(2,2)--(0,2)--(0,0);\end{tikzpicture}}
\newcommand{\subdisconnected}{\begin{tikzpicture}[scale = 0.08]\draw (0,0)--(2,0)--(2,2)--(0,2)--(0,0);\end{tikzpicture}}

\begin{document}

\author{\footnotesize\textbf{Henk Don};
TU Delft, EWI (DIAM), Section Applied Probability,\\ \footnotesize Mekelweg 4, 2628 CD Delft, the Netherlands; e-mail: henkdon@gmail.com}

\title{New methods to bound the critical probability in fractal percolation}
\maketitle

\textbf{Abstract:} We study the critical probability $p_c(M)$ in
two-dimensional $M$-adic fractal percolation. To find lower
bounds, we compare fractal percolation with site percolation.
Fundamentally new is the construction of an computable increasing
sequence that converges to $p_c(M)$. We prove that $p_c(2)> 0.881$
and $p_c(3)>0.784$.

For the upper bounds, we introduce an iterative
random process on a finite alphabet $\mathscr{A}$, which is easier
to analyze than the original process.  We show that $p_c(2)<0.993$, $p_c(3)<0.940$ and $p_c(4)<0.972$.

\section{Introduction}

\begin{figure}[!h]
\begin{tabular}{lllllll}
\includegraphics*[width =2.1cm]{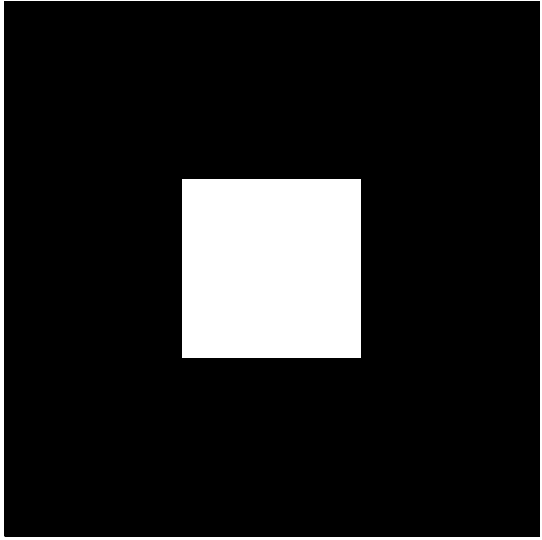}
&
\includegraphics*[width =2.1cm]{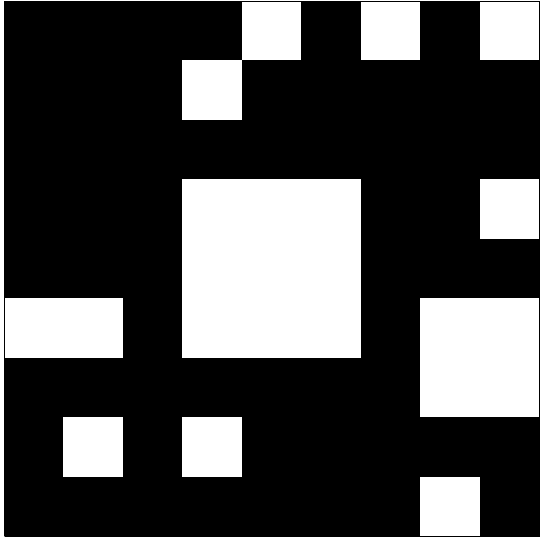}
&
\includegraphics*[width =2.1cm]{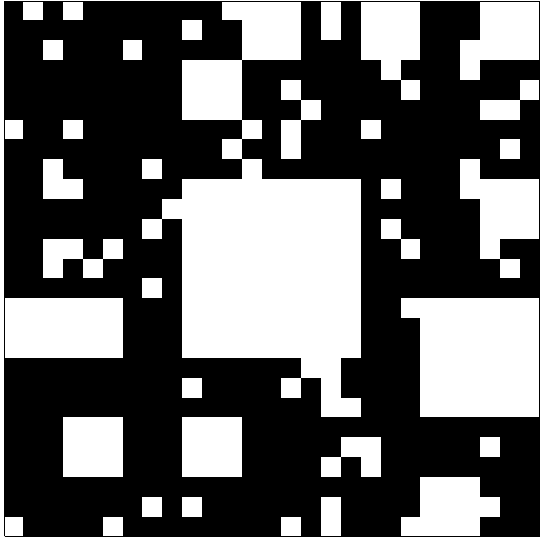}
&
\includegraphics*[width =2.1cm]{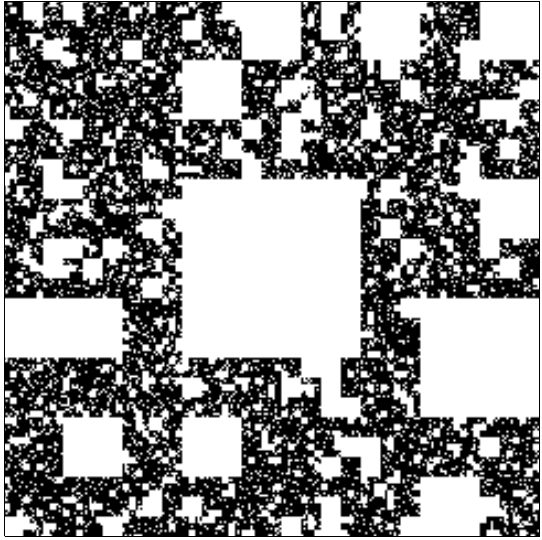}
&
\includegraphics*[width =2.1cm]{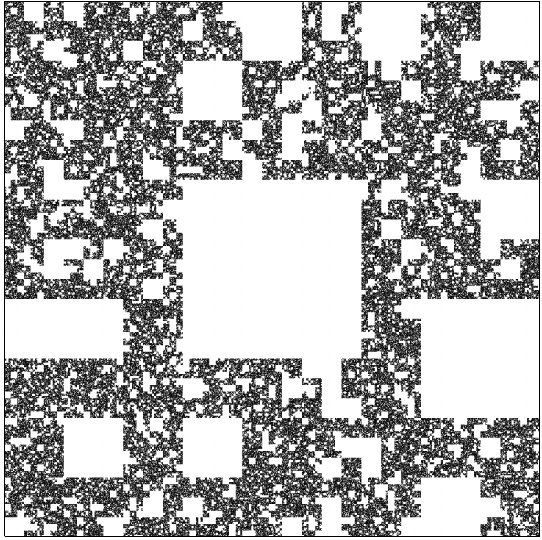}
\end{tabular}\caption{Realizations of $K_n, n=1,2,3,5,7$ for $M=3$ and $p=0.85$.}\label{fig_realizations}
\end{figure}

Fractal percolation has been introduced by Mandelbrot in 1974 as a
model for turbulence and is discussed in his book \emph{The
Fractal Geometry of Nature} \cite{Mandelbrot}. Several equivalent
formal definitions of this process can be found in the literature
(see e.g. \cite{Chayes,DekkingGrimmett,Falconer}). Here we only
give an informal definition of the two-dimensional case. Let $K_0$
be the unit square and choose an integer $M\geq 2$ and a parameter
$p\in [0,1]$. To obtain $K_1$, divide $K_0$ into $M^2$ equal
subsquares, each of which survives with probability $p$ and is
discarded with probability $1-p$, independently of all other
subsquares. Now do the same procedure in all surviving squares, in
order to obtain $K_2$. Iterating this process gives a decreasing
sequence of sets $(K_n)_{n\in\mathbb{N}}$, see Figure
\ref{fig_realizations}. Let $K = \bigcap_{n\in\mathbb{N}}K_n$ be
the limit set.

It was shown in 1988 by Chayes, Chayes and Durrett \cite{Chayes}
that there exists a non-trivial critical value $p_c(M)$ such that
a.s. the largest connected component in $K$ is a point for
$p<p_c(M)$ and with positive probability there is a connected
component intersecting
opposite sides of the unit square for $p\geq p_c(M)$.\\
For all $M\geq 2$, the value of $p_c(M)$ is unknown. Several
attempts have been made to find bounds for $p_c(M)$. It is easy to
see that $K$ is empty a.s. if $p\leq 1/M^2$, which implies $p_c(M)
> 1/M^2$. The argument in \cite{Chayes} is already a bit smarter:
any left-right crossing has to cross the line
$\left\{1/M\right\}\times [0,1]$ somewhere. A crossing of this
line in $K_n$ means that there is a pair of adjacent squares on
opposite sides of this line. Such pairs form a branching process
with mean offspring $p^2M$ and consequently $p_c(M)> 1/\sqrt{M}$.
For the case $M=2$ this was sharpened by White in 2001 to
$p_c(2)\geq 0.810$, who used a set that dominates
$K$ and has a simpler structure to study.\\
Sharp upper bounds are harder to obtain. The first idea to get
rigorous upper bounds for $M\geq 2$ was given by Chayes, Chayes
and Durrett \cite{Chayes}, but (in their own words) these bounds
are ridiculously close to $1$. For $M=3$, they show that
$p_c(3)<0.9999$ (although in fact one can prove that
$p_c(3)<0.993$ with their method), which was improved by Dekking
and Meester \cite{Dekking} to $p_c(3)<0.991$. Chayes et al only
treat $M=3$, but they point out that the same idea works for any
$M\geq 3$. The case $M=2$ can be treated by comparing with $M=4$.
As is noted by van der Wal \cite{Wal}, a coupling argument gives
$p_c(2)\leq 1-(1-\sqrt{p_c(4)}\ )^4$. Following this approach
gives $p_c(4)<0.998$ and $p_c(2)<1-10^{-12}$.\\
In this paper we present ideas to find significantly sharper lower
and upper bounds. To find lower bounds, we compare fractal
percolation with site percolation. In particular, we prove the
following result (we will present a precise definition of
$\pi_n(p,M)$ in section \ref{section_coupling}):

\begin{theorem}\label{theorem_coupling}
Let $M$ be fixed. Define $\pi_n(p) = \mathbb{P}(\textrm{two sides
are connected in }K_n(p))$. If $\pi_n(p) < p_c^{site}$ for some
$n$, then $p < p_c(M)$.
\end{theorem}

This theorem leads to the construction of a increasing computable
sequence $(p_c^n(M))_{n=0}^\infty$ of lower bounds for $p_c(M)$.
However, these computations are quite demanding: to find
$p_c^n(M)$, one needs to consider all possible realizations of
$K_n$. In section \ref{section_classifying} we develop methods to
bound $(p_c^n(M))$ from below by classifying realizations of $K_n$
into some different types, where the set of types $\mathscr{A}$
does not depend on $n$. The fractal percolation iteration process
now induces an iterative random process on $\mathscr{A}$, which is
easier to analyze than the original process. Specifically, the
recursive structure allows us to investigate the limit for large
$n$. Similar ideas are discussed for the upper bounds, but here we
do not need the coupling with site percolation. For the cases
$M=2,3$ and $4$, we use these insights to give computer aided
proofs for the following bounds:
\begin{theorem}
The following bounds hold for $p_c(M), M=2,3,4$:
\begin{enumerate}
\item $p_c(2)>0.881$ and $p_c(3)>0.784$;
\item $p_c(2)<0.993$,
$p_c(3)<0.940$ and $p_c(4)<0.972$.
\end{enumerate}
\end{theorem}

\section{Lower bounds for $p_c(M)$}

In this section we develop methods to calculate lower bounds for
the critical value of two-dimensional fractal percolation. First
we briefly introduce site percolation and then we prove a coupling
with fractal percolation that allows us to find lower bounds for
$p_c(M)$. In particular, we construct an increasing sequence of
lower bounds and we prove that this sequence converges to
$p_c(M)$. At the end of this section we show how to use these
ideas to obtain numerical results.

\subsection{Site percolation}

Consider the infinite two-dimensional square lattice in which each
vertex is open with probability $p$ and closed otherwise. In this
model the percolation probability $\zeta(p)$ is defined as the
probability that the origin belongs to an infinite open cluster.
The critical probability is given by
$$
p^{site}_c := \inf\left\{p:\zeta(p)>0\right\}.
$$
It has been shown by van den Berg and Ermakov \cite{Berg} that
$p_c^{site}>0.556$. The following classical property (see e.g.
\cite{Grimmett} and the references therein) will be used to couple
site percolation to fractal percolation.
\begin{property}\label{property_site}
Take a box of $n\times n$ vertices. Suppose $p < p_c^{site}$. Then
the probability that there is an open cluster intersecting
opposite sides of the box converges to $0$ as
$n\rightarrow\infty$.
\end{property}

\subsection{Coupling site percolation and fractal
percolation}\label{section_coupling}

In the fractal percolation model one usually adopts the following
definitions: a set in the unit square is said to percolate if it
contains a connected component intersecting both the left side and
the right side of the square. Let
$$
\theta_n(p,M) = \mathbb{P}(K_n(p,M)\textrm{\
percolates}),\quad\quad\theta(p,M) = \mathbb{P}(K(p,M)\textrm{\
percolates}).
$$
\piccaption{A percolating set\label{figure_diagonal}}
\parpic[r]{\includegraphics*[width=2.5cm]{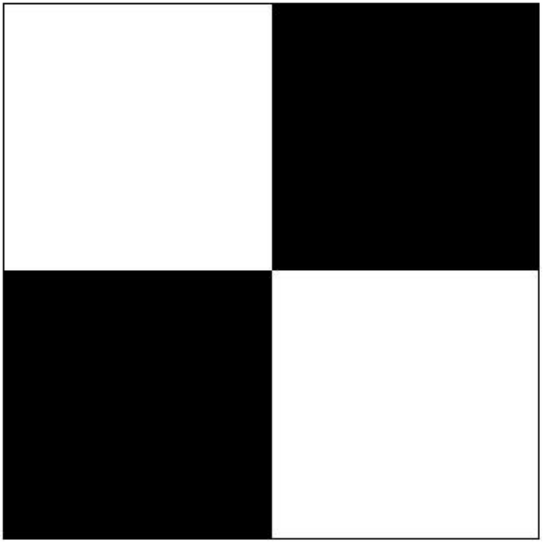}}
The critical probability is defined as
$$
p_c(M) := \inf\left\{p:\theta(p,M)>0\right\}.
$$
We will often suppress some of the dependence on $M$ and $p$. It
is well known (see \cite{Orzechowski}) that
$$
\lim_{n\rightarrow\infty} \theta_n(p) = \theta(p) =
\mathbb{P}(\bigcap_{n=0}^\infty \left\{K_n(p) \textrm{\
percolates}\right\}).
$$
To obtain a proper coupling, we will slightly modify the above
definitions. For example, the set $[0,1/2]^2\cup[1/2,1]^2$
percolates (see Figure \ref{figure_diagonal}). We would like to
ignore such diagonal connections in fractal percolation, since in
site percolation diagonal connections do not exist. Therefore we
redefine percolation as follows:
\begin{definition}\label{def_percolation}
We say $K_n$ percolates if it contains $M^n$-adic squares
$S_1,\ldots,S_k$ such that
\begin{itemize}
\item $S_1$ intersects the left side of the square
$\left\{0\right\}\times[0,1]$.

\item $S_k$ intersects the right side of the square
$\left\{1\right\}\times[0,1]$.

\item $S_i$ shares a full edge with $S_{i+1}$ for $1\leq i\leq
k-1$.
\end{itemize}
\end{definition}

A diagonal connection in $K_n$ can only be present in $K_{n+1}$ if
both subsquares in the corners survive. This means that the
connection is preserved with probability at most $p^2$. As a
consequence, if $K_n$ does not percolate in the sense of
Definition \ref{def_percolation}, then there will be no
percolation in the limit since all (countably many) diagonal
connections break down almost surely. It follows that this
modification of the definition does not change the limiting
percolation probability. From now on also connections, crossings
and connected components in $K_n$ are similarly redefined.
\begin{definition}
Denote the four sides of the unit square by $B_1,\ldots,B_4$. If
$K_n$ contains $M^n$-adic squares $S_1,\ldots,S_k$ such that
\begin{itemize}
\item $S_1$ intersects $B_i$ and $S_k$ intersects $B_j$ for some
$i\neq j$,

\item $S_i$ shares a full edge with $S_{i+1}$ for $1\leq i\leq
k-1$,
\end{itemize}
then we say that two sides are connected in $K_n$.
\end{definition}
%
%
\textbf{Proof of Theorem \ref{theorem_coupling}} First define
\emph{delayed} fractal percolation: $F_{m,n}$ is constructed in
the same way as $K_{m+n}$, the only difference being that we do
not discard any squares in the first $m$ construction steps. So we
first divide the unit square into $M^m\times M^m$ subsquares and
only then we start the fractal percolation process in each of
these squares. Delayed fractal percolation stochastically
dominates fractal percolation.
\begin{figure}[!h]
\begin{tabular}{lll}
\includegraphics*[width =3.7cm]{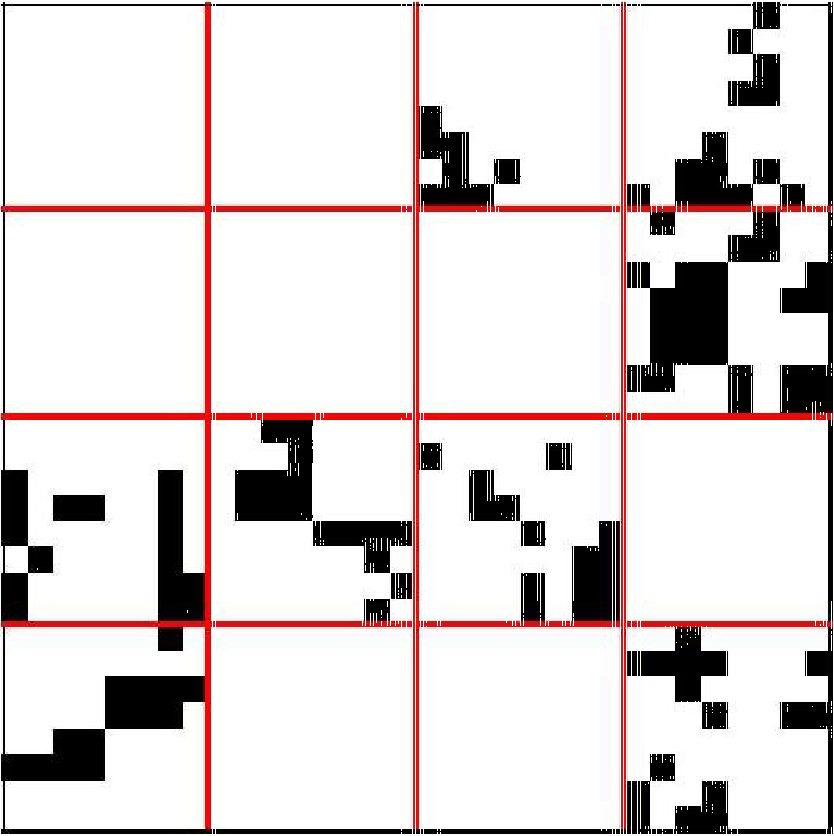}
&
\includegraphics*[width =3.7cm]{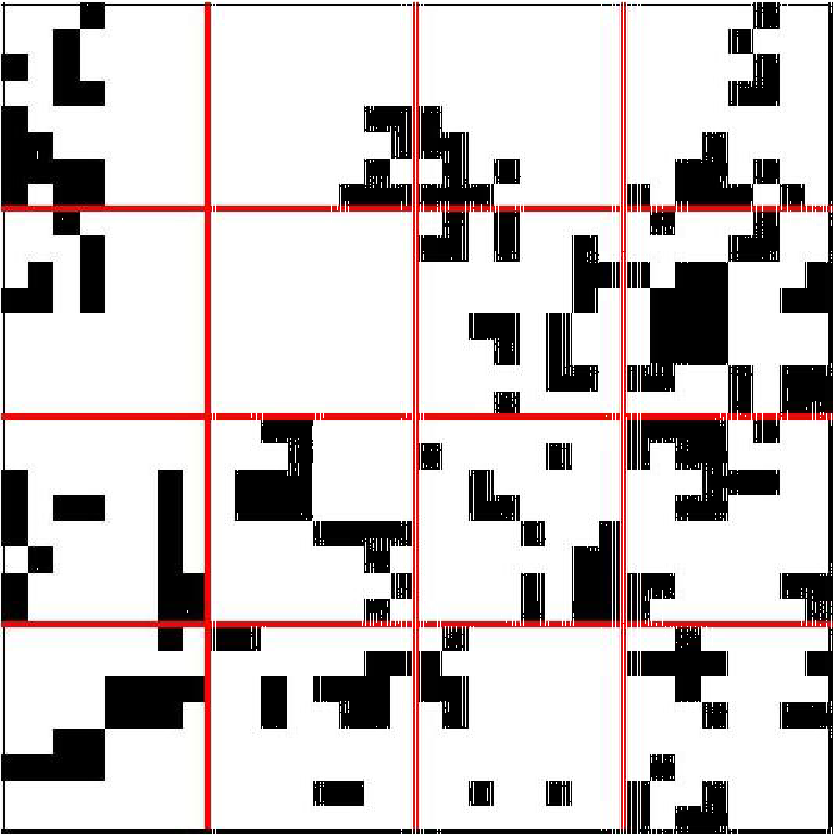}
&
\includegraphics*[width =3.7cm]{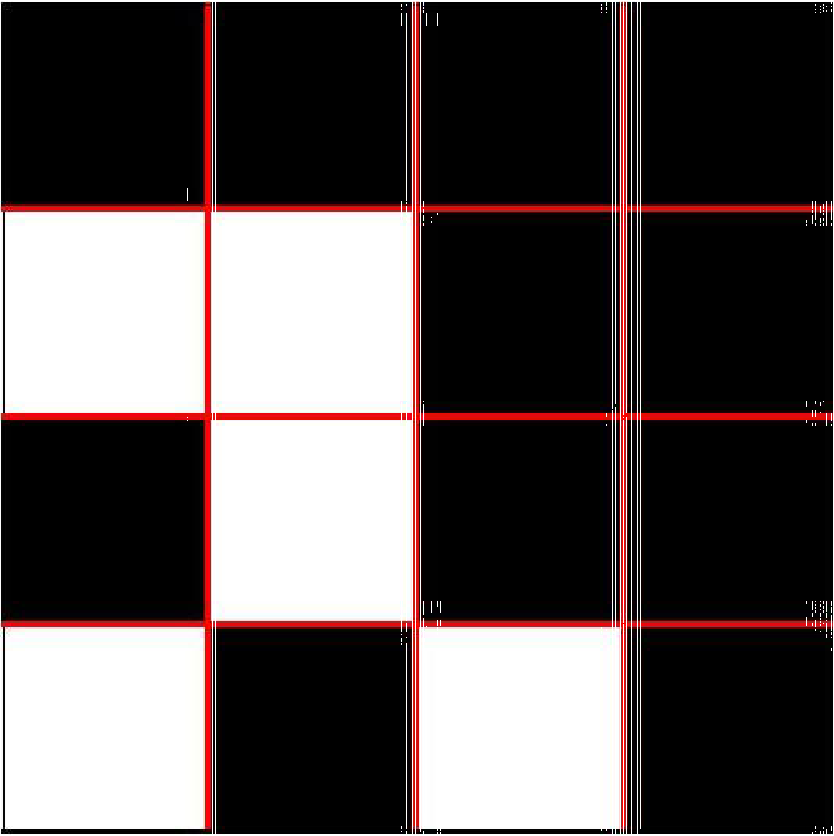}
\end{tabular}\caption{Three coupled realizations to illustrate (\ref{inequality}). Left: $K_5$ for $M=2$ and $p=2/3$. Middle: $F_{2,3}$ for $M=2$ and $p=2/3$. Right: $K_1$ for $M=4$ and $p = \pi_3(2/3,2)$.}\label{fig_inequalities}
\end{figure}
Then we have the following inequalities (illustrated for $M=m=2$
and $n=3$ in Figure \ref{fig_inequalities}):
\begin{equation}\label{inequality}
\theta_{m+n}(p,M) \leq \mathbb{P}(F_{m,n}(p,M) \textrm{
percolates}) \leq \theta_1(\pi_n(p),M^m).
\end{equation}
The first inequality follows from the fact that $F_{m,n}$ is
stochastically larger than $K_{m+n}$. The second inequality can be
explained as follows. Suppose we have two types of squares:
realizations of $K_n(p,M)$ in which two sides are connected (type
1) and realizations in which no sides are connected (type 2).
Suppose we tile a larger square with $M^{2m}$ independent
realizations of $K_n(p,M)$. So the probability on a type 1 square
is $\pi_n(p)$. This larger square is a (scaled) realization of
$F_{m,n}(p,M)$. Now replace all type 1 squares by a full square
and discard all type 2 squares. This gives a realization of
$K_1(\pi(p),M^m)$. Moreover, we claim that this replacement
procedure can not destroy percolation.

To prove this claim, suppose we have a left-right crossing in
$F_{m,n}(p,M)$. This crossing successively traverses $n$th level
squares $S_1,\ldots,S_k$, where $S_i$ shares an edge with
$S_{i+1}$ for $1\leq i\leq k-1$. If $S_i$ is of type 1, then $S_i$
will be replaced by a full square. If $S_i$ is of type 2, then the
crossing enters and leaves $S_i$ at the same side of $S_i$, i.e.
$S_{i-1}=S_{i+1}$ (note that $S_1$ and $S_k$ can not be of type
2). We conclude that if we remove the type 2 squares from
$S_1,\ldots,S_k$, then still each of these squares shares an edge
with its predecessor and successor. From this observation the
claim follows.
%

A first level fractal percolation set can be seen as site
percolation in a finite box. Suppose $\pi_n(p) < p_c^{site}$ for
some $n$ and let the box size $M^m$ tend to $\infty$. By Property
\ref{property_site} we arrive at
\begin{equation}
\lim_{m\rightarrow\infty} \theta_1(\pi_n(p),M^m) = 0.
\end{equation}
Therefore, if $\pi_n(p) < p_c^{site}$, by (\ref{inequality}) we
find that
$$
\lim_{m\rightarrow\infty} \theta_{m+n}(p,M) = 0,
$$
which is equivalent to $\theta(p,M)=0$ and hence $p <
p_c(M)$.\hfill $\Box$

\subsection{A convergent sequence of lower bounds for $p_c(M)$}

In this section we define a sequence of lower bounds for $p_c(M)$.
We prove that this sequence converges to the $p_c(M)$. Let
\begin{equation}\label{eq_sequence}
p_c^n(M) = \sup\left\{p:\pi_n(p,M)<p_c^{site}\right\}.
\end{equation}
Note that $\pi_n(p_c^n)= p_c^{site}$ for all $n$, since $\pi_n(p)$
is continuous in $p$. Since $\pi_n(p)$ is strictly decreasing in
$n$, it follows that $\pi_{n+1}(p_c^n) < p_c^{site}$ for all $n$
and hence by Theorem \ref{theorem_coupling} indeed
$p_c^n(M)<p_c(M)$ for all $n$. The strict monotonicity in $n$ also
implies that $(p_c^n(M))_{n=0}^\infty$ is increasing. The obvious
question now is whether $(p_c^n(M))_{n=0}^\infty$ converges to
$p_c(M)$. We will show that this is indeed the case. First we need
that $\pi_n(p)$ goes to zero if the fractal percolation is
subcritical. Basically this is known, since $K$ is almost surely
disconnected when $p<p_c(M)$ (see \cite{Broman, Chayes}).
\begin{lemma}\label{lemma_convergence}
If $p<p_c(M)$, then $\lim_{n\rightarrow\infty} \pi_n(p,M) = 0$.
\end{lemma}
\textbf{Proof} Suppose $p<p_c(M)$, so $\theta(p)=0$. Note that
a.s. there is an $n$ such that in $K_n$ the two squares in the top
left and bottom left corner are discarded already. Conditioned on
this event, a connection in the limiting set from the left side to
any other side can only occur if it horizontally crosses the
vertical strip $S$ consisting of the squares
$[0,M^{-n}]\times[jM^{-n},(j+1)M^{-1}]$ for $j=1,\ldots,M^n-2$.
From selfsimilarity and subcriticality it follows that in the
limit each of these squares has zero probability to contain a
component connecting opposite sides (horizontally and vertically).
So in $K$ a horizontal crossing of $S$ can only occur if it
crosses a block of two vertically adjacent squares. But as Dekking
and Meester showed (Lemma 5.1 in \cite{Dekking}), $\theta(p)=0$
implies that such a block crossing has zero probability as well.
It follows that connections from the left side to any other side
have zero probability to occur. If this holds for the left side,
then it holds for all sides, and so $\lim_{n\rightarrow\infty}
\pi_n(p) = 0$.\hfill $\Box$\\

The previous lemma makes it easy to prove convergence:

\begin{proposition}\label{proposition_convergence}
The sequence $(p_c^n(M))_{n=0}^\infty$ converges to $p_c(M)$ if
$n\rightarrow\infty$.
\end{proposition}

\textbf{Proof} Let $\varepsilon >0$ and suppose $p =
p_c(M)-\varepsilon$. Then by Lemma \ref{lemma_convergence}, there
exists an $N$ such that $\pi_n(p) < p_c^{site}$ for all $n\geq N$.
Therefore, $p_c^n(M) > p$ for all $n\geq N$. \hfill $\Box$\\

The theory developed so far gives in principle an algorithmic way
to calculate an increasing and converging sequence of lower bounds
for $p_c(M)$: compute the polynomial $\pi_n(p)$ and solve
$\pi_n(p) = p_c^{site}$. Since $p_c^{site}$ is not known exactly,
we replace it in our calculations by the lower bound of van den
Berg and Ermakov. This leads to lower bounds $\tilde p_c^n$ for
$p_c$ that are smaller than $p_c^n$, but they still converge to
$p_c$ (the proof of Proposition \ref{proposition_convergence}
still works).
\begin{example}
\emph{For $M=2$ and $M=3$ we find
\begin{eqnarray*}
\pi_1(p,2) &=& 1 - (1-p)^4,\\
\pi_1(p,3) &=& 1+(1-p)^4(p^5+4p^4(1-p)+6p^3(1-p)^2-1).
\end{eqnarray*}
Solving $\pi_1(p,2) = 0.556$ and $\pi_1(p,3) = 0.556$ gives
$\tilde p_c^1(2)\approx 0.183$ and $\tilde p_c^1(3)\approx 0.178$.
Therefore, $p_c(2)>0.183$ en $p_c(3)>0.178$.}\hfill $\blacksquare$
\end{example}
To find sharper bounds than in the above example, we should take
larger values of $n$. However, for large $n$, the functions
$\pi_n(p)$ are very complicated polynomials, and it is not clear
how to find them in reasonable time. In the next section we will
discuss a way to avoid this problem.

\section{Classifying realizations}\label{section_classifying}

The number of possible realizations of $K_n$ and their complexity
rapidly increases as $n$ grows. In this section we introduce a way
to reduce the complexity without losing too much essential
information on the connectivity structure in $K_n$. Basically, we
will divide the boundary of the unit square in some segments and
the presence or absence of connections between these segments will
determine the type of a realization of $K_n$. We take a finite set
of symbols (also called letters), each representing a type, that does not depend on $n$. This set will be called the
alphabet $\mathscr{A}$. We will analyze the probabilities that
$K_n$ is of a certain type. These ideas can be used to obtain both
lower and upper bounds for $p_c(M)$.

Let $\mathscr{K}_n$ be the set of all possible realizations of
$K_n$. For each $n$ we will define a map $\mathscr{C}_n:\mathscr{K}_n\rightarrow\mathscr{A}$. The sequence of maps $\mathscr{C}=(\mathscr{C}_n)_{n=0}^\infty$ will be called a \emph{classification}. In Section \ref{section_alphabet} we will give a detailed description of the alphabet and classification. For now we only state that $\mathscr{A}$ will be a partially ordered set, having a unique minimum and maximum, with the property that
\begin{equation}\label{eq_conditions}
\mathscr{C}_n(K_n) = \left\{\begin{array}{ll} \min(\mathscr{A}) &
\quad \textrm{if } K_n = \emptyset, \\ \max(\mathscr{A}) &
\quad\textrm{if } K_n = [0,1]^2.
\end{array}\right.
\end{equation}
The letters from $\mathscr{A}$ can be used to create \emph{words}. For our purposes we only need two-dimensional square words of size $M\times M$, denoted by
$$
w = (w_{i,j})_{0\leq i,j\leq M-1} =
\left.
\begin{array}{lll}
w_{0,M-1}&\ldots&w_{M-1,M-1}\\
\vdots&\iddots&\vdots\\
w_{0,0}&\ldots&w_{M-1,0}
\end{array}
\right.,
$$
where $w_{i,j}\in\mathscr{A}$. The set of all such words will be denoted by $\mathscr{A}^{M\times M}$. Since $K_n$ can be obtained by tiling $[0,1]^2$ by scaled realizations of $K_{n-1}$, there is a natural way to associate realizations of $K_n$ to words in $\mathscr{A}^{M\times M}$. First define the \emph{tiles} of $K_n$ as follows
$$
K_n(i,j) := \overline{K_n\cap\left(\left(\frac{i}{M},\frac{i+1}{M}\right)\times\left(\frac{j}{M},\frac{j+1}{M}\right)\right)},\quad 0\leq i,j \leq M-1.
$$
Now rescale and translate them into the unit square:
$$
\hat K_n(i,j) := M\left(K_n(i,j)-\left(\frac{i}{M},\frac{j}{M}\right)\right),\quad 0\leq i,j \leq M-1.
$$
Then each $\hat K_n(i,j)$ either is the empty set, or it can be
seen as a realization of $K_{n-1}$. Now we can map $K_n$ to a word
in $\mathscr{A}^{M\times M}$, provided $\mathscr{C}_{n-1}$ is
known. Define $\mathscr{W}_n:\mathscr{K}_n\rightarrow
\mathscr{A}^{M\times M}$ by
\begin{equation}\label{eq_word}
\mathscr{W}_n(K_n)_{i,j} = \mathscr{C}_{n-1}(\hat K_n(i,j)).
\end{equation}
So far we discussed two maps on $\mathscr{K}_n$: one that maps
realizations to $M\times M$ words (the map $\mathscr{W}_n$) and
one that maps realizations to single letters (the map
$\mathscr{C}_n$). If the word $\mathscr{W}_n(K_n)$ completely
determines $\mathscr{C}(K_n)$, we say the classification is
\emph{regular}:
\begin{definition}
Let $\mathscr{C}=(\mathscr{C}_n)_{n=0}^\infty$ be a
classification. If there exists a map $\phi : \mathscr{A}^{M\times
M}\rightarrow\mathscr{A}$ such that
\begin{equation}\label{eq_regular}
\mathscr{C}_n = \phi\circ\mathscr{W}_n, \quad \quad n\geq 1,
\end{equation}
then we say $\mathscr{C}$ is a \emph{regular classification} and
$\phi$ is called the \emph{word code} of $\mathscr{C}$.
\end{definition}
Note that a regular classification is uniquely defined by its word
code: (\ref{eq_conditions}) defines $\mathscr{C}_0$ and
(\ref{eq_word}) together with the regularity defines
$\mathscr{C}_n$ if $\mathscr{C}_{n-1}$ is known.
\begin{example}\label{example_regular}
\emph{Let $M=2$ and take the alphabet $\mathscr{A} =
\left\{\disconnected,\plus\right\}$, where $\min{\mathscr{A}} =
\disconnected$ and $\max{\mathscr{A}} = \plus$. For
$w\in\mathscr{A}^{2\times 2}$ define
$$
\phi(w) = \left\{
\begin{array}{lll}
\disconnected && \textrm{if } w =
\left.\begin{array}{cc}\disconnected &
\disconnected\\\disconnected & \disconnected\end{array}\right.,\\
\plus && \textrm{otherwise}.
\end{array}
\right.
$$
Let $\mathscr{C}=(\mathscr{C}_n)_{n=0}^\infty$ be regular with
word code $\phi$. Then $\mathscr{C}_0$ is determined by
(\ref{eq_conditions}):
$$
\mathscr{C}_0(K_0(p)) = \mathscr{C}_0([0,1]^2) = \max{\mathscr{A}}
= \plus.
$$
For $n\geq 1, \mathscr{C}_n = \phi\circ\mathscr{W}_n$. For
instance, if $K_1 = [\frac{1}{2},1]^2$, then
\begin{eqnarray*}
\mathscr{C}_1(K_1) &=& \phi(\mathscr{W}_1([\frac{1}{2},1]^2)) =
\phi\left(\begin{array}{cc}\mathscr{C}_0(\hat K_1(1,2)) &
\mathscr{C}_0(\hat K_1(2,2))\\\mathscr{C}_0(\hat K_1(1,1)) &
\mathscr{C}_0(\hat K_1(2,1))\end{array}\right) =\\
&=&\phi\left(\begin{array}{cc}\mathscr{C}_0(\emptyset) &
\mathscr{C}_0([0,1]^2)\\\mathscr{C}_0(\emptyset) &
\mathscr{C}_0(\emptyset)\end{array}\right) =
\phi\left(\begin{array}{cc}\disconnected & \plus\\\disconnected &
\disconnected\end{array}\right) = \plus.
\end{eqnarray*}} \hfill $\blacksquare$
\end{example}

We now want to analyze the probabilities
$\mathbb{P}(\mathscr{C}(K_n(p))=a)$, where $a\in\mathscr{A}$.
Suppose $\mathscr{C}$ is a regular classification with word code
$\phi$. Let
$$
\mathscr{P}_\mathscr{A} = \left\{x\in [0,1]^{|\mathscr{A}|}:
||x||_1 = 1\right\}
$$
be the set of all probability vectors on $\mathscr{A}$. For
$x\in\mathscr{P}_\mathscr{A}$ and $a\in\mathscr{A}$, we denote the
probability that $x$ assigns to $a$ by $x_a$. Take $x\in
\mathscr{P}_\mathscr{A}$, and suppose we construct an $M\times M$
word $w$ in which all letters are chosen independently according
to $x$. Define $F_\mathscr{C}(x) \in \mathscr{P}_\mathscr{A}$ by
$$
(F_\mathscr{C}(x))_a = \mathbb{P}_x(\phi(w)=a),\quad\quad
a\in\mathscr{A}.
$$
The function $F_\mathscr{C}:\mathscr{P}_\mathscr{A}\rightarrow
\mathscr{P}_\mathscr{A}$ will be the key to calculate the
probabilities $\mathbb{P}(\mathscr{C}(K_n(p))=a)$ in an iterative
way, as is shown in the next lemma. Define
$\tau^n(p)\in\mathscr{P}_\mathscr{A}$ by
$$
\tau_a^n(p) = \mathbb{P}(\mathscr{C}(K_n(p))=a).
$$
Let $\tau^{\subdisconnected}$ and $\tau^{\subplus}$ be the vectors
that assign full probability to $\min(\mathscr{A})$ and
$\max(\mathscr{A})$ respectively:
$$
\tau^{\subdisconnected}_a = \left\{
\begin{array}{ll}
1, & a = \min(\mathscr{A})\\
0, & \textrm{otherwise.}
\end{array}
\right.\quad\quad\quad
 \tau^{\subplus}_a = \left\{
\begin{array}{ll}
1, & a = \max(\mathscr{A})\\
0, & \textrm{otherwise.}
\end{array}
\right.
$$
For $S\subseteq \mathscr{A}$, we write
$$
\tau^n_S(p) = \mathbb{P}(\mathscr{C}(K_n(p))\in S) = \sum_{a\in S}
\tau_a^n(p).
$$

\begin{lemma}\label{lemma_recursion}
If the classification $\mathscr{C}$ is regular, then
$$
\tau^{n+1}(p) =
F_\mathscr{C}(p\tau^n(p)+(1-p)\tau^{\subdisconnected}) \quad
\textrm{with initial condition}\quad \tau^0(p) = \tau^{\subplus}.
$$
\end{lemma}
\textbf{Proof} The letters in $\mathscr{W}_{n+1}(K_{n+1}(p))$ are
independent. With probability $p$ a letter corresponds to a scaled
realization of $K_n(p)$, with probability $1-p$ it corresponds to
an empty square. So each letter occurs according to the
probability vector $p\tau^n(p)+(1-p)\tau^{\subdisconnected}$.
Together with the definition of $F_\mathscr{C}$, this gives the
recursion. The initial condition follows from
(\ref{eq_conditions}).\hfill $\Box$\\
\\
This recursion formula is essentially a generalization of the
recursion given in \cite{Dekking}.

\subsection{Strategy for lower bounds}

Recall that our main obstacle in finding sharp bounds for $p_c$ is
that $\pi_n(p)$ is hard to compute. The recursion of Lemma
$\ref{lemma_recursion}$ gives a tool to dominate $\pi_n(p)$ by
something that is easier to compute. The strategy to find lower
bounds for $p_c$ is as follows. Define an alphabet $\mathscr{A}$
with subset $\mathscr{A}_\pi$ and a regular classification
$\mathscr{C}$ (by choosing a word code) in such a way that
$\pi_n(p)\leq\mathbb{P}(\mathscr{C}(K_n(p))\in\mathscr{A}_\pi)$
for all $n$. Now take some $n$ and search for the largest $p$ for
which the latter probability is smaller than $0.556$. Then it
follows that $\pi_n(p)<p_c^{site}$ and therefore $p<p_c$, by
Theorem \ref{theorem_coupling}. We will give an example (which
only gives a very moderate bound) to illustrate this procedure.

\begin{example}\emph{(A lower bound for $M=2$) Let $\mathscr{C}$
be the classification of Example \ref{example_regular}. By
induction it follows that $\mathscr{C}_n(K_n) = \plus$ if $K_n$ is
nonempty. So, if we choose $\mathscr{A}_\pi =
\left\{\plus\right\}$, then definitely
$\pi_n(p)\leq\mathbb{P}(\mathscr{C}(K_n(p))\in\mathscr{A}_\pi) =
\tau^n_{\subplus}(p)$. From the definition of $\phi$ it follows
that in this case
$$
(F_\mathscr{C}(x))_{\subplus}= \mathbb{P}_x(\phi(w)=\plus) =
1-\mathbb{P}_x\left(w =
\begin{array}{cc}
\disconnected & \disconnected\\
\disconnected & \disconnected
\end{array}
\right) = 1-(x_{\subdisconnected})^4 = 1-(1-x_{\subplus})^4.
$$
Note that $(p\tau^n(p)+(1-p)\tau^{\subdisconnected})_{\subplus} =
p\tau_{\subplus}^n(p)$ and apply Lemma \ref{lemma_recursion}:
$$
\tau^{n+1}_{\subplus}(p) =
(F_\mathscr{C}(p\tau^n(p)+(1-p)\tau^{\subdisconnected}))_{\subplus}
= 1-(1-p\tau^{n}_{\subplus}(p))^4,  \quad \textrm{with }
\tau_{\subplus}^0(p) = 1.
$$
Writing $G_p(x):=\tau^1_{\subplus}(px) = 1-(1-px)^4$ leads to
$$
\tau^{n+1}_{\subplus}(p) = G_p(\tau^n_{\subplus}(p)) =
G_p^{n+1}(1).
$$
The function $G_p(x)$ is increasing on $[0,1]$ and $G_p(1)\leq 1$,
so $\tau^n_{\subplus}(p)$ decreases to the largest fixed point of
$G_p$. For $p=0.33$ the fixed point is still below $0.556$ and we
find $\pi_{50}(p) \leq \tau^{50}_{\subplus}(p) \approx 0.554
<0.556 \leq p_c^{site}$. Consequently $p_c(2)>0.33$ by Theorem
\ref{theorem_coupling}.} \hfill $\blacksquare$
\end{example}

\subsection{Strategy for upper bounds}

Our recipe to find upper bounds for $p_c(M)$ is a bit more
involved. We start by defining a partial ordering on the set of
probability vectors $\mathscr{P}_\mathscr{A}$. A set $S\subseteq
\mathscr{A}$ will be called \emph{increasing} if $a\in S$ implies
$b\in S$ for all $b\succeq a$. For
$x,y\in\mathscr{P}_\mathscr{A}$, we now write $x\succeq y$ if
$$
\sum_{a\in S} x_a \geq \sum_{a\in S} y_a, \textrm{for all
increasing } S\subseteq \mathscr{A}.
$$
We say the function $F_\mathscr{C}$ is increasing if
$F_\mathscr{C}(x)\succeq F_\mathscr{C}(y)$ for $x\succeq y$.

\begin{lemma}\label{lemma_taudecreases}
Let $\mathscr{C}$ be a regular classification for which
$F_\mathscr{C}$ is increasing. Then $(\tau^n(p))_{n=0}^\infty$ is
decreasing and $\tau^\infty(p):=
\lim_{n\rightarrow\infty}\tau^n(p)$ exists.
\end{lemma}

\textbf{Proof} Let $S$ be any nonempty increasing subset of
$\mathscr{A}$. Then $\max(\mathscr{A})\in S$, and since
$\tau_{\max(\mathscr{A})}^{\subplus} = 1$ we have
$$
\sum_{a\in S} \tau_a^{\subplus} = 1 \geq \sum_{a\in S} x_a,
\textrm{for all } x\in\mathscr{P}_\mathscr{A}.
$$
Since $\tau^0(p) = \tau^{\subplus}$, it follows that
$$
\tau^0(p) \succeq x, \textrm{for all }
x\in\mathscr{P}_\mathscr{A},
$$
and in particular $\tau^0(p)\succeq \tau^1(p)$. Now we use
induction: suppose $\tau^n(p)\succeq \tau^{n+1}(p)$. Then
$\tau^{n+1}(p)=
F_\mathscr{C}(p\tau^{n}(p)+(1-p)\tau^{\subdisconnected}) \succeq
F_\mathscr{C}(p\tau^{n+1}(p)+(1-p)\tau^{\subdisconnected})=
\tau^{n+2}(p)$ since $F_\mathscr{C}$ is increasing. So
$(\tau^n(p))_{n=0}^\infty$ is decreasing.

To prove existence of $\lim_{n\rightarrow\infty}\tau^n(p)$, we
will show that $\lim_{n\rightarrow\infty}\tau_a^n(p)$ exists for
all $a\in\mathscr{A}$. Let $S_a = \left\{b\in\mathscr{A}:b\succeq
a\right\}$. Then $S_a$ and $S_a\setminus\left\{a\right\}$ are both
increasing sets. Since $(\tau^n(p))_{n=0}^\infty$ is decreasing,
$(\tau_{S_a}^n(p))$ and
$(\tau_{S_a\setminus\left\{a\right\}}^n(p))$ are decreasing
real-valued sequences, bounded from below by $0$. Therefore, their
limits exist and also
$$
\lim_{n\rightarrow\infty} \tau_a^n(p) =
\lim_{n\rightarrow\infty}\left(\tau_{S_a}^n(p)-\tau_{S_a\setminus\left\{a\right\}}^n(p)\right)
$$
exists. Since $\mathscr{A}$ is finite, these limiting
probabilities uniquely determine $\tau^\infty(p)$.\hfill$\Box$\\
\\
To find upper bounds for $p_c(M)$, we want to bound $\theta(p)$
away from $0$. We will construct an alphabet $\mathscr{A}$ with an
increasing subset $\mathscr{A}_\mu\subseteq \mathscr{A}$ and a
regular classification $\mathscr{C}$ for which $F_{\mathscr{C}}$
is increasing. We will do this in such a way that
$\tau_\mu^n(p):=\mathbb{P}(\mathscr{C}(K_n(p))\in\mathscr{A}_\mu)\leq
\theta_n(p)$ for all $n$. If we can prove that $\tau_\mu^\infty(p)
>0$, then it follows that $\theta(p)>0$ and hence $p_c<p$. Finding
$\tau^\infty(p)$ exactly might be not so easy, but the following
lemma gives the key to find a lower bound for it.
\begin{lemma}\label{lemma_upperbound}
Let $\mathscr{C}$ be a regular classification for which
$F_\mathscr{C}$ is increasing.
\begin{enumerate}
\item If $F_\mathscr{C}(px+(1-p)\tau^{\subdisconnected})\succeq x$
for some $x\in\mathscr{P}_\mathscr{A}$ and $p\in(0,1]$, then
$$\tau^\infty(p)\succeq x.$$
\item If in addition $\sum_{a\in\mathscr{A}_\mu}x_a>0$ and
$\theta_n(p)\geq \tau^n_\mu(p)$ for all $n$, then
$$p>p_c.$$
\end{enumerate}
\end{lemma}
\textbf{Proof} We show by induction that $\tau^n(p)\succeq x$ for
all $n\in\mathbb{N}$. First note that $\tau^0(p) = \tau^{\subplus}
\succeq x$. Now suppose $\tau^n(p)\succeq x$ for some $n$. Then by
Lemma \ref{lemma_recursion}
$$
\tau^{n+1}(p) =
F_\mathscr{C}(p\tau^n(p)+(1-p)\tau^{\subdisconnected}) \succeq
F_\mathscr{C}(px+(1-p)\tau^{\subdisconnected}) \succeq x.
$$
Hence indeed $\tau^n(p)\succeq x$ for all $n\in\mathbb{N}$ and
consequently $\tau^\infty(p)\succeq x$.\\
\\
For the second statement, observe that
$$
\theta(p) = \lim_{n\rightarrow\infty}\theta_n(p) \geq
\lim_{n\rightarrow\infty}\tau^n_\mu(p) = \tau^\infty_\mu(p) \geq
\sum_{a\in\mathscr{A}_\mu}x_a>0,
$$
where we used that $\mathscr{A}_\mu$ is an increasing set. Consequently $p_c<p$.\hfill$\Box$\\
\\
The crucial question now is if for given $\mathscr{A}$,
$\mathscr{A}_\mu$ and $\mathscr{C}$ there exists $x\in
\mathscr{P}_\mathscr{A}$ satisfying all requirements and how we
can find it. Here we give a guideline to find numerical results.
First approximate the fixed point $\tau^\infty(p)$ by iterating
the recursion of Lemma \ref{lemma_recursion}. If $p$ is too small,
then $\tau_\mu^\infty(p) = 0$ and no suitable $x$ will exist. If
$p$ is large enough, then $\tau_\mu^\infty(p) > 0$. In the latter
case, if $x$ is an approximation for $\tau^\infty(p)$, we have
$$
\sum_{a\in\mathscr{A}_\mu}x_a > 0\quad\textrm{and}\quad
F_\mathscr{C}(px+(1-p)\tau^{\subdisconnected})\approx x.
$$
Now apply the following trick: let $\varepsilon$ be a small
positive number, then
$$
\sum_{a\in\mathscr{A}_\mu}x_a > 0\quad\textrm{and}\quad
F_\mathscr{C}((p+\varepsilon)x+(1-(p+\varepsilon))\tau^{\subdisconnected})\succeq
x.
$$
Now we have an $x$ that fits into the conditions of Lemma
\ref{lemma_upperbound}, and therefore $p_c<p+\varepsilon$. Before
we give an example of the procedure to find upper bounds for
$p_c(M)$, we will first construct suitable alphabets in the next
section.

\section{Construction of the alphabet and word codes}\label{section_alphabet}

\begin{figure}[!h]
\begin{tabular}{ll}
\includegraphics*[width =6cm]{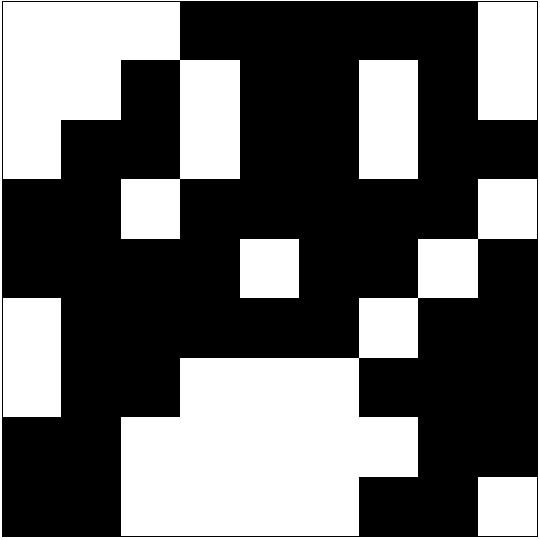}
&
\includegraphics*[width =6cm]{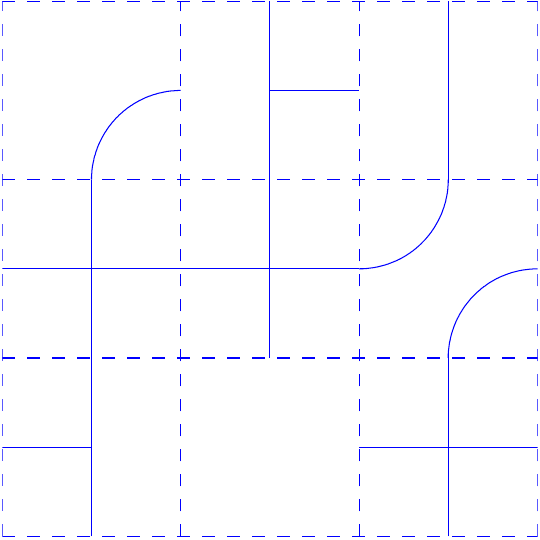}
\end{tabular}
\caption{Left: A realization of $K_2$ for $M=3$. Right: A $3\times
3$ word $w$ over $\mathscr{A}_{3,0}$. If we use the word code
$\Psi_{3,0}$, then this word corresponds to the realization on the
left.}\label{fig_wordplot}
\end{figure}

Let $E=\left\{e_1,\ldots,e_n\right\}$ be a collection of closed
line segments whose union is the boundary of $[0,1]^2$. Assume
that they do not overlap, i.e. the intersection of two segments is
at most a single point. We number them clockwise starting from
$(0,0)$. We will define the alphabet $\mathscr{A}_E$ by means of
non-crossing equivalence relations on $E$.
\begin{definition}
Let $A=\left\{a_1,\ldots,a_n\right\}$ be an ordered set. A
non-crossing equivalence relation on $A$ is a set $R\subseteq
A\times A$ with the following properties:
\begin{enumerate}
\item $(a_i,a_i)\in R$,\quad $i = 1,\ldots,n$, \item $(a_i,a_j)
\in R \Leftrightarrow (a_j,a_i) \in R$,\quad $i,j = 1,\ldots,n$,
\item If $(a_i,a_j),(a_j,a_k)\in R$, then $(a_i,a_k)\in R$,\quad
$i,j,k = 1,\ldots,n$, \item If $i<j<k<l$ and
$(a_i,a_k),(a_j,a_l)\in R$, then $(a_i,a_j)\in R$.
\end{enumerate}
If $(a_i,a_j)\in R$, we say $a_i$ and $a_j$ are equivalent and
write $a_i\sim_R a_j$, or simply $a_i\sim a_j$. The set $[a_i] =
\left\{a_j\in A:a_i\sim a_j\right\}$ is called the equivalence
class of $a_i$.
\end{definition}
The first three properties are the usual reflexivity, symmetry and
transitivity. First consider the simplest case,
$$
E = \left\{e_1,e_2,e_3,e_4\right\} = \left\{\left\{0\right\}\times
[0,1],[0,1]\times\left\{1\right\},\left\{1\right\}\times
[0,1],[0,1]\times\left\{0\right\}\right\}.
$$
An equivalence relation on $E$ can be graphically represented as a
square with some connections (the equivalences) between the
boundaries. For example, the symbol $\trandbl$ represents the
equivalence relation
$$
\left\{(e_1,e_1),(e_1,e_4),(e_2,e_2),(e_2,e_3),(e_3,e_2),(e_3,e_3),(e_4,e_1),(e_4,e_4)\right\},
$$
which has equivalence classes $\left\{e_1,e_4\right\}$ and
$\left\{e_2,e_3\right\}$. Doing this for all non-crossing
equivalence relations on $E$ gives the following set of symbols:
$$
\mathscr{A}_E =
\left\{\disconnected,\topleft,\topright,\bottomright,\bottomleft,\vertical,\horizontal,\tlandbr,\trandbl,\threetop,\threeright,\threebottom,\threeleft,\plus
\right\}.
$$
It will be convenient to identify these symbols with the
corresponding equivalence relations. Now we can easily define a
partial order on $\mathscr{A}_E$: we write $a\preceq b$ if $a$ is
contained in $b$. The alphabet has a unique minimum and maximum,
$\min(\mathscr{A}_E) = \disconnected$ and $\max(\mathscr{A}_E) =
\plus$.\\
\\
In a way similar to the identification of letters with equivalence
relations on $E$, we will identify $M\times M$ words with
equivalence relations on
$$
E_{M\times M} := \bigcup_{0\leq i,j \leq M-1} E+(i,j).
$$
Suppose $w\in\mathscr{A}_E^{M\times M}$. Each letter $w_{i,j}$
from this word is an equivalence relation on $E$, and therefore
$w_{i,j} + (i,j)\times (i,j)$ is an equivalence relation on
$E+(i,j)$. Henceforth,
$$
R_w := \bigcup_{0\leq i,j \leq M-1} w_{i,j} + (i,j)\times (i,j)
$$
is a binary relation on $E_{M\times M}$. In general $R_w$ is not
an equivalence relation. It can be easily checked that $R_w$ is
both reflexive and symmetric. However, $R_w$ might fail to be
transitive. For example, if $w_{0,0} = w_{0,1} = \plus$, then
$$
([0,1]\times\left\{0\right\},[0,1]\times\left\{1\right\})\in R_w
\textrm{\ and\ }
([0,1]\times\left\{1\right\},[0,1]\times\left\{2\right\})\in R_w,
$$
but
$([0,1]\times\left\{0\right\},[0,1]\times\left\{2\right\})\not\in
R_w$. Let $\bar R_w$ be the transitive closure of $R_w$ (that is,
$\overline{R}_w$ is the smallest equivalence relation on
$E_{M\times M}$ that contains $R_w$). Figure \ref{fig_wordplot},
right panel, shows a $3\times 3$ word over $\mathscr{A}_E$. The
sides of the dashed squares are elements of the set $E_{M\times
M}$ and two elements are equivalent with respect to
$\overline{R}_w$ if there is a solid connection between them.

The theory above does not change essentially if we choose $E$ to
be a larger set of line segments that all have equal length. A
special case arises if all segments have length $M^{-n}$. In this
case $E = \left\{e_1,\ldots,e_{4M^n}\right\}$ with $e_i =
\left\{0\right\}\times [(i-1)M^{-n},iM^{-n}]$ for $1\leq i\leq
M^n$ et cetera. We will denote the corresponding alphabet by
$\mathscr{A}_{M,n}$. As a final remark we note that the size of
$\mathscr{A}_E$ is given by the Catalan number
$\frac{1}{|E|+1}{2|E| \choose |E|}$.

\subsection{Weak and strong connectivity}

In this section we will define our word codes. Suppose $w$ is an
$M\times M$ word over the alphabet $\mathscr{A}_{M,n}$, so
$w\in\mathscr{A}_{M,n}^{M\times M}$. Let $E$ be the corresponding
segment set. We define the boundary segment set of $E_{M\times M}$
by
$$
\partial E_{M\times M} := \left\{e\in E_{M\times M}:e\cap\partial
[0,M]^2=e\right\}.
$$
The set $\partial E_{M\times M}$ contains $4M^{n+1}$ segments and
we number them clockwise starting from $(0,0)$:
$e_1^w,\ldots,e_{4M^{n+1}}^w$. Now we partition $\partial
E_{M\times M}$ into the subsets
\begin{equation}\label{eq_partition}
E_i^w = \left\{e^w_{M(i-1)+j}:1\leq j\leq M\right\},\quad 1\leq
i\leq 4M^n.
\end{equation}
Let $\sim$ denote equivalence with respect to $\overline{R}_w$. If
$e\sim f$ for some $e\in E^w_i$ and $f\in E^w_j$, we say $E^w_i$
and $E^w_j$ are \emph{connected}. If $E^w_{i_1},\ldots,E^w_{i_m}$
form a chain of pairwise connected sets, we say $E^w_{i_1}$ and
$E^w_{i_m}$ are \emph{weakly connected}. If there exists an
equivalence class $C$ such that $|C\cap E_i^w|>M/2$ and $|C\cap
E_j^w|>M/2$ or if $i=j$, we say $E^w_i$ and $E^w_j$ are
\emph{strongly connected}. Observe that weak and strong
connectivity are non-crossing equivalence relations on
$\left\{E^w_1,\ldots,E^w_{4M^n}\right\}$.\\
\\
These notions provide tools to define word codes: a word
determines a non-crossing equivalence relation on
$\left\{E^w_1,\ldots,E^w_{4M^n}\right\}$ that will be mapped in
the obvious way to a non-crossing equivalence relation on
$E=\left\{e_1,\ldots,e_{4M^n}\right\}$, which is just a letter in
$\mathscr{A}_{M,n}$. In this way, define word codes
$$
\Phi_{M,n}:\mathscr{A}_{M,n}^{M\times
M}\rightarrow\mathscr{A}_{M,n}\quad\textrm{and}\quad\Psi_{M,n}:\mathscr{A}_{M,n}^{M\times
M}\rightarrow\mathscr{A}_{M,n}
$$
based on weak and strong connectivity respectively. If
$\overline{R}_{w_1} \subseteq \overline{R}_{w_2}$, then
$$
\Phi_{M,n}(w_1)\preceq \Phi_{M,n}(w_2) \quad\textrm{and}\quad
\Psi_{M,n}(w_1)\preceq \Psi_{M,n}(w_2).
$$
Consequently, if $\mathscr{C}$ is a regular classification defined
by one of these word codes, then $F_{\mathscr{C}}$ is increasing.

\begin{example}
\emph{ Consider the word $w$ over $\mathscr{A}_{3,0}$ as shown in
Figure \ref{fig_wordplot}. In this case $\partial E_{M\times M}$
contains $12$ boundary segments, $e_1^w,\ldots,e^w_{12}$. Take
partition sets $E_1^w,\ldots,E_4^w$ as in (\ref{eq_partition}).
For instance,
$$
E_3^w = \left\{e^w_7,e^w_8,e^w_9\right\} =
\left\{\left\{3\right\}\times [2,3],\left\{3\right\}\times
[1,2],\left\{3\right\}\times [0,1]\right\}
$$
contains the three segments at the right side. Then $E_w^1$,
$E_w^2$ and $E_w^4$ are all pairwise connected and $E_w^3$ is
connected to $E_w^4$. Consequently $E_w^i$ is weakly connected to
$E_w^j$ for all $i$
and $j$. Therefore $\Phi_{3,0}(w) = \plus$.\\
Since $e_1^w,e_2^w,e_5^w$ and $e_6^w$ are in the same equivalence
class, $E_1^w$ and $E_2^w$ are strongly connected. There are no
other $i$ and $j$, $i\not =j$ for which $E_i^w$ and $E_j^w$ are
strongly connected. Therefore $\Psi_{3,0}(w)=\topleft$.} \hfill
$\blacksquare$
\end{example}

The idea behind the definitions of $\Phi_{M,n}$ and $\Psi_{M,n}$
is that they guarantee the following key properties:
\begin{property}\label{property_weak} Define a classification
$\mathscr{C}$ by the word code $\Phi_{M,k}$. Then the following
implication holds: if $e,f\in E$ are connected in $K_n$, then
$e\sim f$ in $\mathscr{C}_n(K_n)$.
\end{property}
\textbf{Proof} For $n=0$, we have $\mathscr{C}_0(K_0) =
\max(\mathscr{A}_{M,k})$, which means that all segments are
equivalent. So the statement holds for $n=0$.\\
Now suppose the statement is true for some $n$, and take a
realization of $K_{n+1}$. Let $w = \mathscr{W}_{n+1}(K_{n+1})$. In
each of the tiles $K_{n+1}(i,j)$, the induction hypothesis
applies. So if $e,f\in E_{M\times M}$ are connected in (a scaled,
translated version of) the tile $K_{n+1}(i,j)$, then the
corresponding segments in the letter $w_{i,j}$ are equivalent. If
$e_a,e_b\in E$ are connected in $K_{n+1}$, then this connection
successively traverses some tiles. Hence, there exists segments
$s_1,\ldots,s_m\in E_{M\times M}$ such that $s_1\in E_a^w,s_m\in
E_b^w$ and for which $s_i\sim_w s_{i+1}$, $i=1,\ldots,m-1$ by the
induction hypothesis. Therefore $s_1\sim_w s_m$ and $E_a^w$ is
weakly connected to $E_b^w$. Consequently $e_a\sim e_b$ in
$\Phi_{M,k}(w) = \mathscr{C}_{n+1}(K_{n+1})$.\hfill $\Box$

\begin{property}\label{property_strong}
Define a classification $\mathscr{C}$ by the word code
$\Psi_{M,k}$. Then the reversed implication holds: if $e\sim f$ in
$\mathscr{C}_n(K_n)$, then $e$ and $f$ are connected in $K_n$.
\end{property}

Before proving this property, we introduce some terminology. Let
$C_0 = I$ be an interval and fix integers $M$ and $k>M/2$.
Construct $C_1$ by subdividing $I$ into $M$ subintervals of equal
length and let $k$ of them survive. Repeat this process in each of
the surviving subintervals. The sets in the resulting sequence
$C_0,C_1,\ldots$ will be called $M$-adic fractal majority subsets
of $I$. If $A$ and $B$ are $M$-adic fractal majority subsets of
$I$, then $A\cap B$ contains an interval. If $e,f\in E$ and there
is a connected component in $K_n(M)$ containing $M$-adic fractal
majority subsets of both $e$ and $f$, we say there is a fractal
majority connection between $e$ and $f$.\\
\\
\textbf{Proof of Property \ref{property_strong}} We will actually
prove a stronger statement. Since $K_0 = [0,1]^2$ all boundary
segments are connected to each other in $K_0$ by a fractal
majority connection. Induction hypothesis: if $e\sim f$ in
$\mathscr{C}_n(K_n)$, then $e$ and $f$ are connected
in $K_n$ by a fractal majority connection.\\
Suppose $e_a,e_b\in E$ and $e_a\sim e_b$ in
$\mathscr{C}_{n+1}(K_{n+1})$. Then $E_a^w$ and $E_b^w$ are
strongly connected in $w=\mathscr{W}_{n+1}(K_{n+1})$. With respect
to $\overline{R}_w$ there exists an equivalence class $C\subseteq
E_{M\times M}$ such that $|C\cap E_a^w|>M/2$ and $|C\cap
E_b^w|>M/2$. $C$ contains segments $s_1,\ldots,s_m$ such that
$s_i$ and $s_{i+1}$ are in the same tile of the scaled set
$MK_{n+1}$ for $i = 1,\ldots,m-1$. By the induction hypothesis
there is a fractal majority connection between $s_i$ and $s_{i+1}$
in the corresponding tiles. Since the intersection of two fractal
majority sets is non-empty, all segments $s_1,\ldots,s_m$ are in
the same connected component. Therefore $e_a$ and $e_b$ are
connected in $K_{n+1}$. This connection is a fractal majority
connection since $|C\cap E_a^w|>M/2$ and $|C\cap E_b^w|>M/2$.
\hfill$\Box$\\
\\
Consider the alphabet $\mathscr{A}_{M,k}$ and let $E$ be the
corresponding segment set. Partition $E$ into four sets, each
corresponding to one of the sides of the unit square:
$$
E_i = \left\{e_{1+i|E|/4},\ldots,e_{(1+i)|E|/4}\right\},\quad i =
0,1,2,3.
$$
and define
\begin{equation}\label{eq_percletters}
\begin{array}{lll}
\mathscr{A}_\pi &=& \left\{a\in\mathscr{A}_{M,k}: \exists e\in
E_i, f\in E_j, i\not= j, \textrm{ such that } e\sim_a f
\right\},\\
\mathscr{A}_\mu &=& \left\{a\in\mathscr{A}_{M,k}: \exists e\in
E_1, f\in E_3, \textrm{ such that } e\sim_a f
\right\}.\\
\end{array}
\end{equation}
The following lemma shows that the alphabets and word codes as
defined in this section are suitable for our purposes, see the
discussion in Section \ref{section_classifying}.
\begin{lemma}\label{lemma_strong} Take the alphabet
$\mathscr{A} = \mathscr{A}_{M,k}$ and define $\mathscr{A}_\pi$ and
$\mathscr{A}_\mu$ as in (\ref{eq_percletters}).
\begin{enumerate}
\item Define a classification $\mathscr{C}$ by the word code
$\Phi_{M,k}$. Then $$\tau_{\pi}^n(p)
:=\mathbb{P}(\mathscr{C}(K_n(p))\in\mathscr{A}_\pi)\geq \pi_n(p),
\quad \textrm{ for all }n.$$ \item Define a classification
$\mathscr{C}$ by the word code $\Psi_{M,k}$. Then
$$\tau_{\mu}^n(p)
:=\mathbb{P}(\mathscr{C}(K_n(p))\in\mathscr{A}_\mu) \leq
\theta_n(p),\quad \textrm{ for all }n.$$
\end{enumerate}
\end{lemma}
\textbf{Proof} These statements follow from Property
\ref{property_weak} and \ref{property_strong} respectively.\hfill $\Box$\\
\\
Now we are ready to give an example illustrating how to find upper
bounds. We will keep the example as simple as possible, so that it
can be checked by hand. Therefore our alphabet will contain only
two letters and we will use a simplification of the word code
$\Psi_{3,0}$. Nevertheless, it leads to a bound that already
improves upon the best bound known so far.

\begin{example} \emph{(An upper bound for $M=3$) Let $\mathscr{A} = \left\{\min(\mathscr{A}_{3,0}),\max(\mathscr{A}_{3,0})\right\} =
\left\{\disconnected,\plus\right\}$ and let $E$ be the boundary
segment set corresponding to $\mathscr{A}_{3,0}$. Partition
$\partial E_{M\times M}$ into four sets as in
(\ref{eq_partition}). Define a regular classification
$\mathscr{C}$ by the word code
$$
\phi(w) = \left\{
\begin{array}{lll}
\plus && \textrm{if } E_i^w \textrm{ and } E_j^w \textrm{ strongly connected for all } i,j,\\
\disconnected && \textrm{otherwise}.
\end{array}
\right.
$$
This classification is increasing. Define $\mathscr{A}_\mu =
\left\{\plus\right\}$, which is an increasing subset of
$\mathscr{A}$. Since $\phi(w)\subseteq \Psi_{3,0}(w)$ for all
$w\in \mathscr{A}^{M\times M}$, the second statement of Lemma
\ref{lemma_strong} also applies to this classification. For
$x\in\mathscr{P}_{\mathscr{A}}$ given by
$p\tau^{\subplus}+(1-p)\tau^{\subdisconnected}$, write $x =
(p,1-p)$. Counting all words for which $E_w^1,\ldots,E_w^4$ are
strongly connected gives
$$
(F_\mathscr{C}(x))_{\subplus} =
\mathbb{P}_x(\phi(w)=\plus)=p^9+9p^8(1-p)+20p^7(1-p)^2.
$$
Now choose $p=0.984$ and $y=(0.9720,0.028)$. Then
$(F_\mathscr{C}(y))_{\subplus}\approx 0.9721 > y_{\subplus}$ and
hence $F_\mathscr{C}(y)\succeq y$. Since
$\sum_{a\in\mathscr{A}_\mu}y_a = y_{\subplus}>0$ all conditions
for Lemma \ref{lemma_upperbound} are fulfilled. Hence
$p_c(3)<0.984$.} \hfill $\blacksquare$
\end{example}

\subsection{Monotonicity and convergence}

So far we developed some tools to find bounds for $p_c(M)$. One
would expect that taking larger alphabets results in sharper
bounds, since we can approximate the connectivity structure in
$K_n$ more accurately. In this section we show that this is indeed
the case and that the lower bounds even convergence to $p_c(M)$ if
the alphabet size goes to infinity.\\
For the word code $\Phi_{M,k}$ over $\mathscr{A}_{M,k}$, define
the corresponding classification and let $\tau_\pi^n(p)$ be
defined as before. Then define a critical value as follows:
$$
p_c(\Phi_{M,k}) := \sup\left\{p:\tau_\pi^\infty(p)<p_c^{site}\right\}.
$$
Let $\mathscr{A} = \mathscr{A}_{M,k}$ and define $\mathscr{C}$ by
the word code $\Psi_{M,k}$. Let $\tau_\mu^n(p)$ be as before. Also
here we define a critical value:
$$
p_c(\Psi_{M,k}) := \inf\left\{p:\tau^\infty_\mu(p)>0\right\}.
$$
Now we have the following proposition:
\begin{proposition}
The sequence $(p_c(\Phi_{M,k}))_{k=0}^\infty$ is increasing and
$(p_c(\Psi_{M,k}))_{k=0}^\infty$ is decreasing. Moreover,
$$
\lim_{k\rightarrow\infty} p_c(\Phi_{M,k}) = p_c(M).
$$
\end{proposition}
\textbf{Proof} Denote the segment set corresponding to the
alphabet $\mathscr{A}_{M,k}$ by $E_{M,k} =
\left\{e_1,\ldots,e_{4M^k}\right\}$. Let $\partial
E_{M,k}^{M\times M}$ be the set of boundary segments of $M\times
M$ words over $\mathscr{A}_{M,k}$. Partition $\partial
E_{M,k}^{M\times M}$ into subsets $E_1,\ldots,E_{4M^k}$ as in
(\ref{eq_partition}). Each of these partition sets contains $M$
segments, and the union of all $4M^{k+1}$ segments is equal to
$\partial[0,M]^2$. There is a one-to-one correspondence between
these segments and the elements of $E_{M,k+1}$: if $e\in
E_{M,k+1}$, then $Me\in\partial E_{M,k}^{M\times M}$. Define
$$
C_i = \left\{e\in E_{M,k+1}:Me\in E_i^w\right\}, \quad i =
1,\ldots, 4M^k.
$$
For $e_i\in E_{M,k}$, we have $e_i = \bigcup_{e\in C_i} e$. The
elements of $C_i$ will be called the children of their parent
$e_i$. Let $\mathscr{A}_\pi^{k}$ and
$\mathscr{A}^k_\mu$ be defined according to (\ref{eq_percletters}).\\
Define regular classifications $\mathscr{C}^k$ and
$\mathscr{C}^{k+1}$ by the word codes $\Phi_{M,k}$ and
$\Phi_{M,k+1}$. Denote the corresponding probability vectors by
$^k\tau$ and $^{k+1}\tau$. By induction on $n$ it follows that if
two segments are equivalent in $\mathscr{C}_n^{k+1}(K_n)$, then
their parents are equivalent in $\mathscr{C}_n^{k}(K_n)$.
Therefore, if $\mathscr{C}_n^{k+1}(K_n)\in\mathscr{A}_\pi^{k+1}$
then $\mathscr{C}_n^{k}(K_n)\in\mathscr{A}_\pi^k$. Consequently,
${^{k}}\tau_\pi^n(p)\geq {^{k+1}}\tau_\pi^n(p)$ and so
$$
p_c(\Phi(M,k+1))\geq p_c(\Phi(M,k)).
$$
Now define $\mathscr{C}^k$ and $\mathscr{C}^{k+1}$ by the word
codes $\Psi_{M,k}$ and $\Psi_{M,k+1}$. By induction on $n$: if
$e_i,e_j\in E_{M,k}$ are equivalent in $\mathscr{C}_n^{k}(K_n)$,
then the sets of children $C_i$ and $C_j$ are strongly connected
in $\mathscr{C}_n^{k+1}(K_n)$. So, if
$\mathscr{C}_n^{k}(K_n)\in\mathscr{A}_\mu^{k}$ then
$\mathscr{C}_n^{k+1}(K_n)\in\mathscr{A}_\mu^{k+1}$. Henceforth,
${^k}\tau_\mu^n(p)\leq {^{k+1}}\tau_\mu^n(p)$, so
$$
p_c(\Psi(M,k+1))\leq p_c(\Psi(M,k)).
$$
Having shown the monotonicity of the two sequences, we now turn to
the convergence of $p_c(\Psi(M,k)$. Take the alphabet
$\mathscr{A}_{M,k}$ and define $\mathscr{C}^k$ by the word code
$\Phi_{M,k}$. A realization of $K_n$ consists of $M^n\times M^n$
squares, and letters in $\mathscr{A}_{M,k}$ have $M^k$ boundary
segments at each side. This means that for $n\leq k$ the
classification describes the connectivity structure exactly: two
segments in $\mathscr{C}^k(K_n)$ are equivalent if and only if
they are connected in $K_n$. So $\pi_n(p) = {^k}\tau_\pi^n(p)$ if
$n\leq k$. This implies that for $n\leq k$ we can rewrite
(\ref{eq_sequence}):
\begin{eqnarray*}
p_c^n(M) &=& \sup\left\{p:{^k}\tau_\pi^n(p)<p_c^{site}\right\}\\
&\leq & \sup\left\{p:{^k}\tau_\pi^\infty(p)<p_c^{site}\right\} =
p_c(\Phi_{M,k}) \leq p_c(M),
\end{eqnarray*}
where we used that ${^k}\tau_\pi^n(p)$ decreases in $n$ by Lemma
\ref{lemma_taudecreases}. Proposition
\ref{proposition_convergence} states that $p_c^n(M)$ converges to
$p_c(M)$, so we conclude that $\lim_{k\rightarrow\infty}
p_c(\Phi_{M,k}) = p_c(M)$.\hfill $\Box$

\section{Numerical results}

In this section we present our numerical results. The recursion of
Lemma \ref{lemma_recursion} is the main tool to perform the
calculations. Our implementation in Matlab (everything available
from the author on request) gives the following results:

\begin{proposition}\label{theorem_weakbounds}
Take the alphabet $\mathscr{A}_{M,k}$ and define a classification
$\mathscr{C}$ by the word code $\Phi_{M,k}$. Let $n=1000$ and
define $\tau_{\pi}^n(p)$ as before. Then
\begin{itemize}
\item For $M=2$ and $k=0$, we have $\tau_{\pi}^n(0.785) <
p_c^{site}$. \item For $M=2$ and $k=1$, we have
$\tau_{\pi}^n(0.859) < p_c^{site}$. \item For $M=3$ and $k=0$, we
have $\tau_{\pi}^n(0.715) < p_c^{site}$.
\end{itemize}
\end{proposition}

\begin{corollary}\label{corollary_weakbounds}
$p_c(2)>0.859$ and $p_c(3)>0.715$.
\end{corollary}

\begin{figure}[!h]
\begin{tabular}{l}
\includegraphics*[width =10cm]{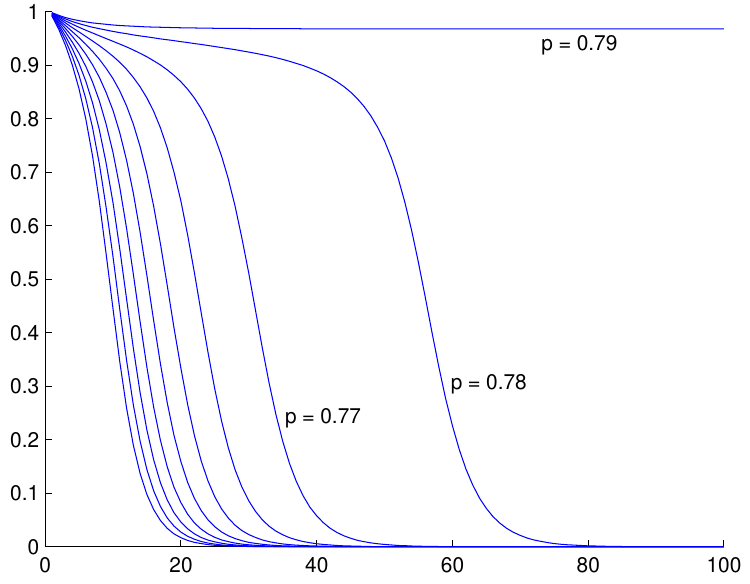}
\end{tabular}\caption{Plot of $\tau_\pi^n(p)$ for $p=0.7+0.01k$ where $k=0,\ldots, 9$ as functions of $n$. Especially note the difference between $p=0.78$ and $p=0.79$.}\label{fig_plottau}
\end{figure}

\textbf{Proof} This follows from Lemma \ref{lemma_strong} and
Theorem \ref{theorem_coupling}.\hfill $\Box$\\
\\
Figure \ref{fig_plottau} illustrates for the case $M=2$ and $k=0$
how $\tau_{\pi}^n(p)$ behaves as a function of $n$ for some values
of $p$. The values of $\tau_{\pi}^n(p)$ were calculated by
iterating the recursion of Lemma \ref{lemma_recursion}.

For larger values of $k$ the computations were too complicated to
perform in a reasonable computation time. For example, the segment
set $E_{2,2}$ contains $16$ segments and therefore the alphabet
$\mathscr{A}_{2,2}$ already contains $\frac{1}{17}{32\choose 16}=
35357670$ letters. Nevertheless we will explain that it is
possible to improve the bounds of Corollary
\ref{corollary_weakbounds} by taking other alphabets or word codes.\\
For $M=2$, define the segment set $E$ by dividing the left and
right side of $[0,1]^2$ into four segments of lenght $1/4$ and the
bottom and top side into two segments of length $1/2$. This leads
to an alphabet $\mathscr{A}_E$ that is in some sense in between
$\mathscr{A}_{2,1}$ and $\mathscr{A}_{2,2}$ and will be denoted by
$\mathscr{A}_{2,3/2}$. Analogous to our previous approach, we
define a classification by choosing the word code that is based on
weak connectivity. For this classification we find
$\tau_{\pi}^{50}(0.876)<p_c^{site}$, which implies
$p_c(2)>0.876$.\\
One can improve this even a bit more by taking the alphabet
$\mathscr{A}_{2,2}$ and defining a word code $\tilde\Phi_{2,2}$
that is a bit simpler than $\Phi_{2,2}$ as follows. If at least
one of the letters in a $2\times 2$ word $w$ equals
$\min(\mathscr{A}_{2,2})$, then $\tilde\Phi_{2,2}(w) =
\Phi_{2,2}(w)$. Otherwise, define $\tilde\Phi_{2,2}(w)$ by first
mapping each of the four letters to $\mathscr{A}_{2,3/2}$ and then
mapping the new word to $\mathscr{A}_{2,2}$, in both steps using
weak connectivity. This simplifies the required calculations a
lot, and leads to $\tilde\tau_\pi^{200}(0.881)<p_c^{site}$. Since
$\tilde\Phi_{2,2}(w)\supseteq \Phi_{2,2}(w)$, we have
$\tilde\tau_\pi^n(p) \geq \tau_\pi^n(p) \geq \pi_n(p)$ for all
$n$. We conclude that $p_c(2)>0.881$.\\
For $M=3$ we improved the lower bound of Corollary
\label{corollary_weakbounds} by using segments of length $1/3$ at
the left and the right side of $[0,1]^2$ and segments of length
$1$ at the bottom and the top side. Denote the resulting alphabet
by $\mathscr{A}_{3,1/2}$ and choose the word code based on weak
connectivity. This gives $\tau_{\pi}^{100}(0.784)<p_c^{site}$,
whence $p_c(3)>0.784$.\\
These calculations have been checked by Arthur Bik, a mathematics
student at Delft University of Technology. He independently
implemented the algorithms and reproduced all results, except the
bound $p_c(2)>0.881$. This was due to the fact that his program
was not fast enough to perform the calculations in a reasonable
time. Concluding, the best lower bounds we found are

\begin{theorem}
$p_c(2)>0.881$ and $p_c(3)>0.784$.
\end{theorem}

Now let us turn to the upper bounds. The strategy described in
Section \ref{section_classifying} leads to the following results:

\begin{proposition}\label{theorem_strongbounds}
Take the alphabet $\mathscr{A}_{M,k}$ and define a classification
$\mathscr{C}$ by the word code $\Psi_{M,k}$. Let $n=1000$. The
conditions
$F_{\mathscr{C}}(px+(1-p)\tau^{\subdisconnected})\succeq x$ and
$\sum_{a\in\mathscr{A}_\mu}x_a>0$ hold if $x$ and $p$ are chosen
as follows:
\begin{itemize}
\item For $M=3$ and $k=0$, choose $p = 0.958$ and $x =
\tau^n(0.9579)$. \item For $M=4$ and $k=0$, choose $p = 0.972$ and
$x = \tau^n(0.9719)$.
\end{itemize}
\end{proposition}

\begin{corollary}\label{corollary_strongbounds}
$p_c(3)<0.958$ and $p_c(4)<0.972$.
\end{corollary}

\textbf{Proof} This follows from Lemma \ref{lemma_upperbound} and
Lemma \ref{lemma_strong}.\hfill $\Box$\\
\\
For $M=3$, the result can be sharpened by using the alphabet
$\mathscr{A}_{3,1/2}$. The classification is again defined by
strong connectivity. In that case the choice $p = 0.940$ and $x =
\tau^{1000}(0.9399)$
satisfies all conditions, so $p_c(3)<0.940$.\\
The algorithm for $M=4$ can be slightly adapted to find a bound
for $M=2$. Each realization of $K_n$ for $M=4$ can be seen as a
realization of $K_{2n}$ for $M=2$. Therefore, we can still use the
word code $\Psi_{4,0}$. The only thing that changes is the way the
probabilities are computed. Given $\tau^n(p)$, the letters in the
$4\times 4$ word $w=\mathscr{W}_n(K_{n+2})$ occur according to the
following rule: The word $w$ consists of four $2\times 2$ blocks.
In each of these blocks either all letters are equal to
$\min(\mathscr{A})$ (with probability $1-p$) or they are
independent of each other chosen according to
$p\tau^n(p)+(1-p)\tau^{\subdisconnected}$ (with probability $p$).
Basically we collapse two construction steps of $K_n$ into one
step. These ingredients determine the recursion. Performing the
calculations we find that the conditions are satisfied for $p =
0.993$ and $x = \tau^{1000}(0.9929)$,
henceforth $p_c(2)<0.993$.\\
Also for the upper bounds Arthur Bik checked our results. He
independently reproduced our bounds, except for the bound
$p_c(3)<0.940$ (for similar reasons as before). Summarizing, our
best upper bounds are
\begin{theorem}
$p_c(2)<0.993$, $p_c(3)<0.940$ and $p_c(4)<0.972$.
\end{theorem}

\textbf{Acknowledgement} The author thanks Michel Dekking for the
entertaining discussions on this topic and for his useful remarks
on this paper. I am also grateful to Arthur Bik, who has spent a
lot of time checking the numerical results.


\begin{thebibliography}{3}

\bibitem{Berg} Berg, J. van den and Ermakov, A. -- \emph{A New Lower Bound for the Critical Probability of Site Percolation on the Square
Lattice}, Random Structures and Algorithms, Vol. 8, No. 3 (1996).

\bibitem{Broman} Broman, E. and Camia, F. -- \emph{Universal Behavior of Connectivity Properties in Fractal Percolation
Models}, Electronic Journal of Probability, 15 (2010), 1394--1414.

\bibitem{Chayes} Chayes, J.T., Chayes, L. and Durrett, R. -- \emph{Connectivity Properties of Mandelbrot's Percolation Process}, Probability Theory and Related Fields, 77 (1988), 307--324.

\bibitem{DekkingGrimmett} Dekking, F.M. and Grimmett, G.R. -- \emph{Superbranching processes and projections of random Cantor
sets}, Probability Theory and Related Fields, 78 (1988), 335--355.

\bibitem{Dekking} Dekking, F.M. and Meester, R.W.J. -- \emph{On the Structure of Mandelbrot's Percolation and Other Random Cantor Sets}, Journal of Statistical Physics, Vol. 58, Nos. 5/6 (1990), 1109--1126.

\bibitem{Falconer} Falconer, K. -- \emph{Random fractals}, Math.
Proc. Cambr. Phil. Soc. 100 (1986), 559-582.

\bibitem{Grimmett} Grimmett G.R. -- \emph{Percolation}, Second edition, Springer-Verlag, Berlin (1999).

\bibitem{Mandelbrot} Mandelbrot, B.B. -- \emph{The Fractal Geometry of
Nature}, Freeman, San Francisco, 1983.

\bibitem{Orzechowski} Orzechowski, M. -- \emph{Geometrical and topological properties of fractal
percolation}, PhD Thesis, University of St. Andrews, 1997.

\bibitem{Wal} Wal, P. van der -- \emph{Random substitutions and fractal percolation}, PhD thesis, Delft University of Technology, 2002.

\bibitem{White} White, D.G. -- \emph{On the Value of the Critical Point in Fractal Percolation}, Random Structures and Algorithms, 18 (2001), 332--345.
\end{thebibliography}
\end{document}